\documentclass[final,numberedheadings]{aipproc}
\layoutstyle{6x9}
\usepackage{amsfonts}
\usepackage{amsmath}

\newtheorem{theorem}{Theorem}[section]

\newtheorem{definition}[theorem]{Definition}

\begin{document}

\title{A Note on $3$-quasi-Sasakian Geometry}


\author{Beniamino Cappelletti Montano}{
  address={Dipartimento di Matematica,
Universit\`{a} degli Studi di Bari, Via E. Orabona 4, 70125 Bari, Italy \\
cappelletti@dm.uniba.it, dileo@dm.uniba.it}}

\author{Antonio De Nicola}{
  address={CMUC, Department of Mathematics, University of Coimbra, 3001-454 Coimbra, Portugal\\
adenicola@mat.uc.pt}}

\author{Giulia Dileo}{
  address={Dipartimento di Matematica,
Universit\`{a} degli Studi di Bari, Via E. Orabona 4, 70125 Bari, Italy \\
cappelletti@dm.uniba.it, dileo@dm.uniba.it}}

\begin{abstract}
$3$-quasi-Sasakian manifolds were recently studied by the authors
as a suitable setting unifying $3$-Sasakian and $3$-cosymplectic geometries.
In this paper some geometric properties of this class of almost
$3$-contact metric manifolds are briefly reviewed, with an emphasis
on those more related to physical applications.

\end{abstract}

\maketitle

%%%%%%%%%%%%%%%%%%%%%%%%%%%%%%%%%%%%%%%%%%%%
%% MAINMATTER
%%%%%%%%%%%%%%%%%%%%%%%%%%%%%%%%%%%%%%%%%%%%

\section{Introduction}
The class of $3$-quasi-Sasakian manifolds is the analogue in the setting of $3$-structures
of the class of quasi-Sasakian manifolds, introduced by Blair \cite{blair0}
and later studied among others by Tanno \cite{tanno}, Kanemaki \cite{kanemaki1},
Olszak \cite{olszak2}. More recent are the examples of applications of quasi-Sasakian manifolds
to string theory found by Friedrich and his collaborators \cite{agricola1,friedrich}.
Just like quasi-Sasakian manifolds include Sasakian and cosymplectic
manifolds, so $3$-quasi-Sasakian manifolds unify $3$-Sasakian and $3$-cosymplectic geometry.
A $3$-quasi-Sasakian manifold can arise, for example, as the product
of a $3$-Sasakian manifold and a hyper-K\"{a}hler manifold (see Sect.~\ref{ranksection} or \cite{mag}).
The setting of $3$-structures has been recently the object of a wider interest from both
mathematicians and physicists due to the important role acquired by the $3$-Sasakian and the related
quaternionic structures in supergravity and superstring theory, where they appear
in the so called hypermultiplet solutions (see e.~g. \cite{acharya,agricola1,gibbons,yee}). This note
contains a concise review of the main properties of $3$-quasi-Sasakian manifolds, recently studied by
the authors in \cite{mag}, together with some relevant properties of the two important subclasses
of $3$-Sasakian and $3$-cosymplectic manifolds which were compared in \cite{cappellettidenicola}.

\section{$3$-quasi-Sasakian geometry}
An \emph{almost contact metric manifold} is a $(2n+1)$-dimensional
manifold $M$ endowed with a field $\phi$ of endomorphisms of the
tangent spaces, a vector field $\xi$, called  \emph{Reeb vector field},
a $1$-form $\eta$ satisfying $\phi^2=-I+\eta\otimes\xi$, $\eta\left(\xi\right)=1$
(where $I \colon TM\rightarrow TM$ is the identity mapping)
and a \emph{compatible} Riemannian metric $g$ such that $g\left(\phi X,\phi Y\right) =
g\left(X,Y\right)-\eta\left(X\right)\eta\left(Y\right)$ for all
$X,Y\in\Gamma\left(TM\right)$.
%From the definition it follows also that $\phi\xi=0$, $\eta\circ\phi=0$
%and that the $(1,1)$-tensor field $\phi$ has constant rank $2n$(cf. \cite{blair1}).
The manifold is said to be \emph{normal} if the tensor field
$N^{(1)}=[\phi,\phi]+2d\eta\otimes\xi$ vanishes identically.
The $2$-form $\Phi$ on $M$ defined
by $\Phi\left(X,Y\right)=g\left(X,\phi Y\right)$ is called the
\emph{fundamental $2$-form} of the almost contact metric manifold
$\left(M,\phi,\xi,\eta,g\right)$.
    Normal almost contact metric manifolds such that both $\eta$ and $\Phi$ are
closed are called \emph{cosymplectic manifolds} and those such that $d\eta=\Phi$ are called
\emph{Sasakian manifolds}.
    The notion of  quasi-Sasakian structure  unifies those of Sasakian and cosymplectic
structures. A \emph{quasi-Sasakian manifold} is defined as a
normal almost contact metric manifold whose fundamental $2$-form
is closed. A quasi-Sasakian manifold $M$ is said to be of rank $2p$ (for some $p\leq n$) if
$\left(d\eta\right)^p\neq 0$ and $\eta\wedge\left(d\eta\right)^p=0$ on $M$,
and to be of rank $2p+1$ if $\eta\wedge\left(d\eta\right)^p\neq 0$ and
$\left(d\eta\right)^{p+1}=0$  on $M$ (cf. \cite{blair0,tanno}).
Blair proved that there are no quasi-Sasakian manifolds of even rank.
Just like Blair and Tanno did, we will only consider quasi-Sasakian manifolds of constant (odd) rank.
If the rank of $M$ is $2p+1$, then the module $\Gamma(TM)$ of vector
fields over $M$ splits into two submodules as follows:
$\Gamma(TM)={\cal E}^{2p+1}\oplus{\cal E}^{2q}$, $p+q=n$, where
\(
{\cal E}^{2q}=\{X\in\Gamma(TM)\; | \; i_X d\eta=0 \mbox{ and } i_X \eta=0\}
\)
and  ${\cal E}^{2p+1}={\cal
E}^{2p}\oplus\left\langle\xi\right\rangle$, ${\cal E}^{2p}$ being
the orthogonal complement of ${\cal
E}^{2q}\oplus\left\langle\xi\right\rangle$ in
$\Gamma\left(TM\right)$. These modules satisfy $\phi {\cal
E}^{2p}={\cal E}^{2p}$ and $\phi {\cal E}^{2q}={\cal E}^{2q}$ (cf.
\cite{tanno}).

We now come to the main topic of our paper, i.e. $3$-quasi-Sasakian geometry,
which is framed into the more general setting of almost $3$-contact geometry.
An  \emph{almost $3$-contact metric manifold}  is a
$\left(4n+3\right)$-dimensional smooth  manifold $M$ endowed with
three almost contact structures $\left(\phi_1,\xi_1,\eta_1\right)$,
$\left(\phi_2,\xi_2,\eta_2\right)$,
$\left(\phi_3,\xi_3,\eta_3\right)$ satisfying the following
relations, for any even permutation
$\left(\alpha,\beta,\gamma\right)$ of $\left\{1,2,3\right\}$,
\begin{gather}
\phi_\gamma=\phi_{\alpha}\phi_{\beta}-\eta_{\beta}\otimes\xi_{\alpha}=-\phi_{\beta}\phi_{\alpha}+\eta_{\alpha}\otimes\xi_{\beta},\\
\nonumber
\xi_{\gamma}=\phi_{\alpha}\xi_{\beta}=-\phi_{\beta}\xi_{\alpha}, \ \
\eta_{\gamma}=\eta_{\alpha}\circ\phi_{\beta}=-\eta_{\beta}\circ\phi_{\alpha},
 \label{3-sasaki}
\end{gather}
and a Riemannian metric $g$ compatible with each of them. It is well
known that in any almost $3$-contact metric manifold the Reeb vector
fields $\xi_1,\xi_2,\xi_3$ are orthonormal with respect to the
compatible metric $g$ and that the structural group of the tangent
bundle is reducible to $Sp\left(n\right)\times I_3$. 
Due to a general result (cf. \cite[Prop.~3.6.2]{joyce}), it follows that
any $3$-quasi-Sasakian manifold is a \emph{spin manifold}, i.~e. there
exists a double cover of the orthonormal frame bundle $SO(TM)$, which is
non-trivial on the fibers of the latter, by a principal bundle with structure group
a spin group. Then, to each representation of the spin group corresponds an associated
vector bundle whose sections are called \emph{spinor fields} (see \cite{lawson} for details).
The existence of spinor fields is  required in quantum theories which encompass fermions.

Moreover, by putting ${\cal{H}}=\bigcap_{\alpha=1}^{3}\ker\left(\eta_\alpha\right)$ one
obtains a $4n$-dimensional \emph{horizontal} distribution on $M$ and the tangent
bundle splits as the orthogonal sum $TM={\cal{H}}\oplus{\cal{V}}$, where
${\cal V}=\left\langle\xi_1,\xi_2,\xi_3\right\rangle$ is the \emph{vertical} distribution.

\begin{definition}
A \emph{$3$-quasi-Sasakian manifold} is an almost
$3$-contact metric manifold $(M,\phi_\alpha,\xi_\alpha,\eta_\alpha,g)$
such that each almost contact structure is quasi-Sasakian.
\end{definition}

The class of $3$-quasi-Sasakian manifolds includes as special cases the well-known
$3$-Sasakian and $3$-cosymplectic manifolds.

The following theorem combines the results obtained in Theorems 3.4 and 4.2 of \cite{mag}.
\begin{theorem}\label{classification}
Let $(M,\phi_\alpha,\xi_\alpha,\eta_\alpha,g)$ be a
$3$-quasi-Sasakian manifold. Then the $3$-dimensional distribution
$\cal V$ generated by $\xi_1$, $\xi_2$, $\xi_3$ is integrable.
Moreover, $\cal V$ defines a totally geodesic and Riemannian
foliation of $M$ and for any even permutation
$(\alpha,\beta,\gamma)$ of $\left\{1,2,3\right\}$ and for some $c\in \mathbb R$
\begin{equation*}
\left[\xi_\alpha,\xi_\beta\right]=c\xi_\gamma.
\end{equation*}
\end{theorem}

Using Theorem \ref{classification} we may divide
$3$-quasi-Sasakian manifolds in two classes according
to the behaviour of the leaves of the
foliation $\cal V$: those $3$-quasi-Sasakian manifolds for which
each leaf of $\cal V$ is locally $SO\left(3\right)$ (or
$SU\left(2\right)$) (which corresponds to take in Theorem
\ref{classification} the constant $c\neq 0$), and those for which
each leaf of $\cal V$ is locally an abelian group (this corresponds
to the case $c=0$).

The preceding theorem also allows to define a canonical metric connection
on any $3$-quasi-Sasakian manifold. Indeed, let $\nabla^B$ be the
Bott connection associated to $\mathcal V$, that is the partial
connection on the normal bundle $TM/{\mathcal V}\cong\mathcal H$ of
$\mathcal V$ defined by
\(
\nabla^B_{V}Z:=[V,Z]_{\mathcal H}
\)
for all $V\in\Gamma(\mathcal V)$ and $Z\in\Gamma(\mathcal H)$.
Following \cite{tondeur} we may construct an adapted connection on
$\mathcal H$ putting
\begin{equation*}
\tilde\nabla_{X}Y:=\left\{
                     \begin{array}{ll}
                       \nabla^B_{X}Y, & \hbox{if $X\in\Gamma(\mathcal V)$;} \\
                       (\nabla_{X}Y)_{\mathcal H}, & \hbox{if $X\in\Gamma(\mathcal H)$.}
                     \end{array}
                   \right.
\end{equation*}
 This connection can be also
extended to a connection on all $TM$ by requiring that
$\tilde\nabla\xi_\alpha=0$ for each $\alpha\in\left\{1,2,3\right\}$.
Some properties of this global connection have been considered in
\cite{cappellettidenicola} for any almost $3$-contact metric
manifold. Now combining Theorem \ref{classification} with \cite[Theorem
3.6]{cappellettidenicola} we have:
\begin{theorem}\label{connessionecanonica}
Let $(M,\phi_\alpha,\xi_\alpha,\eta_\alpha,g)$ be a
$3$-quasi-Sasakian manifold. Then there exists a unique metric connection
$\tilde\nabla$ on $M$ satisfying the following properties:
\begin{enumerate}
  \item[(i)]
  $\tilde\nabla\eta_\alpha=0$, $\tilde\nabla\xi_\alpha=0$, for each
  $\alpha\in\left\{1,2,3\right\}$,
  \item[(ii)]
  $\tilde T\left(X,Y\right)=2\sum_{\alpha=1}^{3}d\eta_{\alpha}(X,Y)\xi_\alpha$,
  for all $X,Y\in\Gamma\left(TM\right)$.
\end{enumerate}
\end{theorem}

\section{The rank of a $3$-quasi-Sasakian manifold}\label{ranksection}
For a $3$-quasi-Sasakian manifold one can consider the ranks of the three  structures $(\phi_\alpha,\xi_\alpha,\eta_\alpha,g)$.
The following theorem assures that these three ranks coincide.
\begin{theorem}[\cite{mag}]\label{rango}
Let $(M,\phi_\alpha,\xi_\alpha,\eta_\alpha,g)$ be a
$3$-quasi-Sasakian manifold of dimension $4n+3$. Then the $1$-forms $\eta_1$, $\eta_2$
and $\eta_3$ have the same rank $4l+3$
or $4l+1$, for some $l\leq n$, according to
$\left[\xi_\alpha,\xi_\beta\right]=c\xi_\gamma$ with $c\neq 0$, or
$\left[\xi_\alpha,\xi_\beta\right]=0$, respectively.
\end{theorem}

According \ to \ Theorem \ref{rango}, \ we \ say \ that \ a \
$3$-quasi-Sasakian \
 manifold $(M,\phi_\alpha,\xi_\alpha,\eta_\alpha,g)$ has \emph{rank}
 $4l+3$ or $4l+1$ if any quasi-Sasakian structure has such rank. We may thus classify $3$-quasi-Sasakian manifolds of dimension
$4n+3$, according to their rank. For any $l\in\{0,\ldots,n\}$ we
have one class of manifolds such that
$\left[\xi_\alpha,\xi_\beta\right]=c\xi_\gamma$ with $c\neq 0$,
and one class of manifolds with
$\left[\xi_\alpha,\xi_\beta\right]=0$. The total number of classes
amounts then to $2n+2$.
In the following we will use the notation ${\cal E}^{4m}:=\{X\in\Gamma({\cal H})\; | \; i_X d\eta_\alpha=0\}$, while
${\cal E}^{4l}$ will be the orthogonal complement of ${\cal E}^{4m}$ in $\Gamma({\cal H})$,
${\cal E}^{4l+3}:={\cal E}^{4l} \oplus \Gamma({\cal V})$, and
${\cal E}^{4m+3}:={\cal E}^{4m} \oplus \Gamma({\cal V})$.

We now consider the class of $3$-quasi-Sasakian manifolds such that
$[\xi_\alpha,\xi_\beta]=c\xi_\gamma$ with $c\neq 0$ and let $4l+3$ be the rank.
In this case, according to \cite{blair0}, we define for each structure
$(\phi_\alpha,\xi_\alpha,\eta_\alpha,g)$ two $(1,1)$-tensor fields
$\psi_\alpha$ and $\theta_\alpha$ by putting
\begin{equation*}
 \psi_\alpha X=\left\{
             \begin{array}{ll}
               \phi_\alpha X, & \hbox{if $X\in {\cal{E}}^{4l+3}$;}\\
               0, & \hbox{if $X\in {\cal{E}}^{4m}$;}
             \end{array}
              \right.
\ \textrm{  } \
 \theta_\alpha X=\left\{
             \begin{array}{ll}
               0, & \hbox{if $X\in {\cal{E}}^{4l+3}$;}\\
               \phi_\alpha X, & \hbox{if $X\in {\cal{E}}^{4m}$.}\\
             \end{array}
              \right.
\end{equation*}
Note that, for each $\alpha\in\left\{1,2,3\right\}$ we have
$\phi_\alpha=\psi_\alpha+\theta_\alpha$.
Next, we define a new (pseudo-Riemannian, in general) metric
$\bar{g}$ on $M $ setting
\begin{equation*}
\bar{g}\left(X,Y\right)=\left\{
                          \begin{array}{ll}
                            -d\eta_\alpha\left(X,\phi_\alpha Y\right), & \hbox{for $X,Y\in{\cal E}^{4l}$;} \\
                            g\left(X,Y\right), & \hbox{elsewhere.}
                          \end{array}
                        \right.
\end{equation*}
This definition is well posed by virtue of normality and of \cite[Lemma 5.3]{mag}.
$(M,\phi_\alpha,\xi_\alpha,\eta_\alpha,\bar{g})$
is in fact a hyper-normal almost $3$-contact metric manifold, in general
non-$3$-quasi-Sasakian.
We are now able to formulate the following decomposition theorem, proven in \cite{mag}.
\begin{theorem}
Let $(M^{4n+3},\phi_\alpha,\xi_\alpha,\eta_\alpha,g)$ be a
$3$-quasi-Sasakian manifold of rank $4l+3$ with $\left[\xi_\alpha,\xi_\beta\right]=2\xi_\gamma$. Assume $[\theta_\alpha, \theta_\alpha]=0$ for some $\alpha\in\{1,2,3\}$ and $\bar{g}$ positive definite on ${\cal E}^{4l}$.
Then $M^{4n+3}$ is locally the  product of a $3$-Sasakian manifold $M^{4l+3}$ and a hyper-K\"{a}hlerian manifold $M^{4m}$ with $m=n-l$.
\end{theorem}

    We now consider the class of $3$-quasi-Sasakian manifolds such that
$[\xi_\alpha,\xi_\beta]=0$ and let $4l+1$ be the rank.
In this case we define for each structure
$(\phi_\alpha,\xi_\alpha,\eta_\alpha,g)$ two $(1,1)$-tensor fields
$\psi_\alpha$ and $\theta_\alpha$ by putting
\begin{equation*}
 \psi_\alpha X=\left\{
             \begin{array}{ll}
               \phi_\alpha X, & \hbox{if $X\in {\cal{E}}^{4l}$;}\\
               0, & \hbox{if $X\in {\cal{E}}^{4m+3}$;}
             \end{array}
              \right.
\ \textrm{  } \
 \theta_\alpha X=\left\{
             \begin{array}{ll}
               0, & \hbox{if $X\in {\cal{E}}^{4l}$;}\\
               \phi_\alpha X, & \hbox{if $X\in {\cal{E}}^{4m+3}$.}\\
             \end{array}
              \right.
\end{equation*}
Note that for each $\alpha$ the maps $-\psi_\alpha^2$ and
$-\theta_\alpha^2+ \eta_\alpha\otimes\xi_\alpha$ define an almost product
structure which is integrable if and only if $[-\psi_\alpha^2,-\psi_\alpha^2]=0$
or, equivalently, $[\psi_\alpha,\psi_\alpha]=0$. Under this assumption the structure
turns out to be $3$-cosymplectic:
\begin{theorem}[\cite{mag}]
Let \ $(M,\phi_\alpha,\xi_\alpha,\eta_\alpha,g)$  be  a
$3$-quasi-Sasakian  manifold of rank $4l+1$ such  that
$\left[\xi_\alpha,\xi_\beta\right]=0$ for any
$\alpha,\beta\in\{1,2,3\}$ and $[\psi_\alpha, \psi_\alpha]=0$ for
some $\alpha\in\{1,2,3\}$. Then $M$ is  a $3$-cosymplectic
manifold.
\end{theorem}

    As we have remarked before, $3$-Sasakian and $3$-cosymplectic manifolds
belong to the class of $3$-quasi-Sasakian manifolds, having respectively rank $4n+3=\dim(M)$ and rank $1$.
We now briefly collect some additional properties of these two important subclasses.
We have seen that the vertical distribution $\cal{V}$ is integrable already in any $3$-quasi-Sasakian manifold.
Ishihara (\cite{ishihara}) has shown that if the foliation defined by $\cal{V}$
is regular then the space of leaves is a quaternionic-K\"{a}hlerian manifold.
Boyer, Galicki and  Mann have proved the following more general result.
\begin{theorem}[\cite{galicki1}]\label{proiezione}
Let $\left(M^{4n+3},\phi_\alpha,\xi_\alpha,\eta_\alpha,g\right)$
be a $3$-Sasakian manifold such that the Killing vector fields
$\xi_1$, $\xi_2$, $\xi_3$ are complete. Then
\begin{description}
    \item[(i)]
    $M^{4n+3}$ is
    an Einstein manifold of positive scalar curvature equal to
    $2\left(2n+1\right)\left(4n+3\right)$.
    \item[(ii)] Each leaf of the foliation $\cal{V}$
    is a $3$-dimensional homogeneous spherical space form.
    \item[(iii)] The space of leaves $M^{4n+3}/\cal{V}$ is a quaternionic-K\"{a}hlerian
    orbifold of dimension $4n$ with positive scalar
    curvature equal to $16n\left(n+2\right)$.
\end{description}
\end{theorem}
We consider now the horizontal distribution: on the one hand, in the
$3$-Sasakian subclass ${\cal{H}}$ is never integrable. On the other
hand, in any $3$-cosymplectic manifold $\cal H$ is integrable since
each $\eta_\alpha$ is closed. Furthermore, the projectability with
respect to $\cal{V}$ is always granted, as the following theorem
shows.
\begin{theorem}[\cite{cappellettidenicola}]\label{hyperkahler}
Every regular $3$-cosymplectic manifold projects onto a hyper-K\"{a}hlerian manifold.
\end{theorem}

As a corollary, it follows that every $3$-cosymplectic manifold is Ricci-flat.

In \cite{cappellettidenicola} the horizontal flatness of such
structures has been studied. In particular it has been proven to
be equivalent to the existence of Darboux-like coordinates, that is
local coordinates $\left\{x_1,\ldots,x_{4n},z_1,z_2,z_3\right\}$ with respect to
which, for each $\alpha\in\left\{1,2,3\right\}$, the fundamental
$2$-forms $\Phi_\alpha=d\eta_\alpha$ have constant components and
$\xi_\alpha=a^1_\alpha \frac{\partial}{\partial z_1} +
a^2_\alpha\frac{\partial}{\partial z_2} + a^3_\alpha
\frac{\partial}{\partial z_3}$,  $a_\alpha^\beta$ being functions
depending only on the coordinates $z_1,z_2,z_3$. Consequently, in view of Theorem
\ref{proiezione} and Theorem \ref{hyperkahler} we have the following result.
\begin{theorem}[\cite{cappellettidenicola}]\label{sasakiano}
A $3$-Sasakian manifold does not admit any Darboux-like coordinate system.
On the other hand, a $3$-cosymplectic manifold admits a Darboux-like coordinate system
around each of its points if and only if it is flat.
\end{theorem}

\section{Final Remarks}
A number of natural questions arose during the development of our work on $3$-quasi-Sasakian manifolds.
We have seen that $3$-Sasakian manifolds do not admit any Darboux coordinate system,
while on $3$-cosymplectic manifolds such coordinate exist if and only if the manifold is flat, so
it is natural to wonder whether these coordinates do not exist on any $3$-quasi-Sasakian manifold of rank greater than one.
Another important topic would be to study the projectability of $3$-quasi-Sasakian manifolds for understanding the general relation
between this class and the quaternionic structures, since the $3$-Sasakian manifolds project over quaternionic-K\"{a}hler
structures while the structure of the leaf space turns out to be globally hyper-K\"{a}hlerian in the $3$-cosymplectic case.
Finally, as both $3$-Sasakian and $3$-cosymplectic manifolds are Einstein manifolds a natural question would be to ask
whether all $3$-quasi-Sasakian manifolds are Einstein. However, since we have already found an example of an $\eta$-Einstein,
non-Einstein $3$-quasi-Sasakian manifold in \cite{mag}, the natural problem now becomes to establish if there is any
$3$-quasi-Sasakian manifolds which is not $\eta$-Einstein. We will try to address some of these questions in the next future.

\section*{Acknowledgments}
The second author acknowledges financial support by a CMUC postdoctoral fellowship.

\end{document}